\documentclass[11pt]{article}
\usepackage{mathrsfs}
\usepackage{epic,eepic,epsf,epsfig}
\usepackage{amsfonts,srcltx,mathrsfs}
\usepackage{amsmath}
\usepackage{amssymb}
\usepackage{amsbsy}
\usepackage{graphicx}
\usepackage{amsfonts}
\usepackage{color}
\usepackage{mathrsfs}
\usepackage{epic,eepic,epsf,epsfig}
\usepackage{graphicx}
\usepackage{amsfonts,srcltx,mathrsfs}
\usepackage{array}
\usepackage{amsthm}
\usepackage{algorithm}
\usepackage{algorithmic}

\newtheorem{theorem}{Theorem}

\newtheorem{clm}{Claim}

\def\f{\noindent}

\begin{document}

\markboth{ et. al}{Antimagic orientations of disconnected even regular graphs}

\title{Antimagic orientations of disconnected even regular graphs
}

\author{Chen Song, Rong-Xia Hao\footnote{Corresponding author. Email: rxhao@bjtu.edu.cn (R.-X. Hao)}\\[0.2cm]
{\em\small Department of Mathematics, Beijing Jiaotong University,}\\ {\small\em
Beijing 100044, P.R. China}}

\date{}
\maketitle

 A $labeling$ of a digraph $D$ with $m$ arcs is a bijection from the set of arcs of $D$ to $\{1,2,\ldots,m\}$. A labeling of $D$ is $antimagic$ if no two vertices in $D$ have the same vertex-sum, where the vertex-sum of a vertex $u \in V(D)$ for a labeling is the sum of labels of all arcs entering $u$ minus the sum of labels of all arcs leaving $u$. An antimagic orientation $D$ of a graph $G$ is $antimagic$ if $D$ has an antimagic labeling. Hefetz, M$\ddot{u}$tze and Schwartz in [J. Graph Theory 64(2010)219-232] raised the question: Does every graph admits an antimagic orientation? It had been proved that for any integer $d$, every 2$d$-regular graph with at most two odd components has an antimagic orientation. In this paper, we consider the 2$d$-regular graph with many odd components. We show that every 2$d$-regular graph with any odd components has an antimagic orientation provide each odd component with enough order.

\medskip

\f {\em Keywords:} Regular graph; Antimagic labeling; Antimagic orientation.

\section{Introduction}
All graphs in this paper are finite and simple. For a graph $G$, let $|G|$ denote the number of vertices of $G$. For a path $P$, let $|P|$ denote the length of $P$. For an orientation $D$ of graph $G$, $D$ is a digraph, we use $A(D)$ and $V(D)$ to denote the set of arcs and vertices of $D$, respectively. We define $[i,j]:=\{i,i+1,\ldots,j\}$, for any two positive integers $i$ and $j$. A $labeling$ of $D$ with $m$ arcs is a bijection from $A(D)$ to $[1,m]$. A labeling of $D$ is $antimagic$ if no two vertices in $D$ have the same vertex-sum, where the vertex-sum of a vertex $u \in V(D)$ for a labeling is the sum of labels of all arcs entering $u$ minus the sum of labels of all arcs leaving $u$. An antimagic orientation $D$ of $G$ is $antimagic$ if $D$ has an antimagic labeling. A graph $G$ has an antimagic orientation if an antimagic orientation of $G$ is antimagic. Let $D$ be an orientation of a graph $G$ with $m$ edges. For any labeling $c: A(D)\rightarrow[1,m]$ of $D$ and any vertex $u \in V(D)$, we use $s_{D}(u)$ to denote the $vertex$-$sum$ of $u$ for the labeling $c$ of orientation $D$.

Hefetz, M$\ddot{u}$tze and Schwartz~\cite{H} raised the question: Does every graph admits an antimagic orientation? For this question and  any integer $d \geq 1$, they proved the following solutions: $(a)$ every $(2d-1)$-regular graph admits an antimagic orientation;
$(b)$ every connected $2d$-regular graph $G$ admits an antimagic orientation if $G$ has a matching covers all but at most one vertex of $G$. Alon et al~\cite{A} obtained that the dense graphs are antimagic. Cranston~\cite{CR} proved that regular bipartite graphs are antimagic.
Chang et al.~\cite{CH} discussed the antimagic labeling of regular graphs. Cranston et al~\cite{CRL} proved that regular graphs of odd degree are antimagic. Recently, Shan et al.~\cite{S} support this conjecture by proving that every biregular bipartite graph admits an antimagic orientation. Li et al.\cite{L} proved that every connected $2d$-regular graph admits an antimagic orientation. Let $G$ be a $2d$-regular graph, where $d \geq 2$ is an integer. The result that $G$ admits an antimagic orientation if $G$ has at most two odd components is proved in~\cite{L}.

It remained unknown whether every disconnected $2d$-regular graph admits an antimagic orientation. In this paper, first we find an orientation of the disconnected $2d$-regular graph, then we find a labeling based on this orientation and finally we show that this labeling is antimagic provide each component with enough order. The main results of this paper are Theorem~\ref{thm1} and Theorem~\ref{thm2}.

\begin{theorem}\label{thm1}
For any integer $d \geq 2$, let $G$ be a $2d$-regular graph with components $G_{1},G_{2},\ldots,\\G_{q}$, where $G_{1},G_{2},\ldots,G_{k}$ are odd components such that $|G_{1}|\leq|G_{2}|\leq \ldots \leq|G_{k}|$. Then the following results hold.

$(1)$ For $k \in [0,5d+4]$, $G$ admits an antimagic orientation.

$(2)$ For $k \geq 5d+5$, if $|G_{1}|\geq 2x_{0}+5$, then $G$ admits an antimagic orientation.
Where $x_{0}$ is the unique positive integer solution for one of equations: $k=(2d-2)(x+2)+0$, $k=(2d-2)(x+2)+1$,\ldots, $k=(2d-2)(x+2)+(d+8), k=(2d-2)(x+1)+(d+9)$, $k=(2d-2)(x+1)+(d+10),\ldots,k=(2d-2)(x+1)+(2d-3)$.
\end{theorem}

From Theorem~\ref{thm1}, the following result is derived directly.

\begin{theorem}\label{thm2}
For any integer $d \geq 2$, let $G$ be a $2d$-regular graph with $q$ components. If each odd component of $G$ has enough order, then
$G$ admits an antimagic orientation.
\end{theorem}

\section{Proof of Theorem~\ref{thm1}}

A closed walk in a graph is an {\it Euler tour} if it traverses every edge of the graph exactly once. There is a Euler theorem that a connected graph admits an Euler tour if and only if every vertex has even degree.

In this section, we study the antimagic orientations of disconnected $2d$-regular graphs for any integer $d\geq 2$. We will prove Theorem~\ref{thm1} by the following process: Firstly, we find an orientation $D^{\ast}$ of the given graph $G$. Secondly, we label the edges of $G$ by three algorithms. Thirdly, we show that the orientation $D^{\ast}$ is antimagic by proving the labeling of $D^{\ast}$ which is given by the algorithm is antimagic.

\begin{figure}[!ht]\label{fig1}
\begin{center}
\includegraphics[scale=0.6]{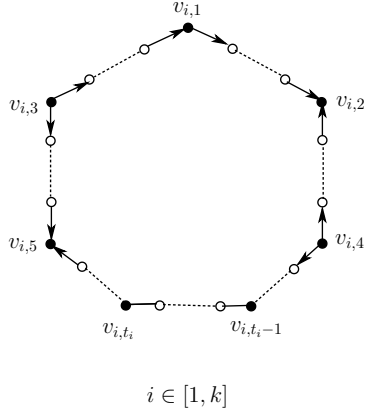}
\end{center}
\vskip-0.2cm
\caption{The orientation of an odd cycle}\label{F1}
\end{figure}

\begin{figure}[!ht]
\begin{center}
\includegraphics[scale=0.6]{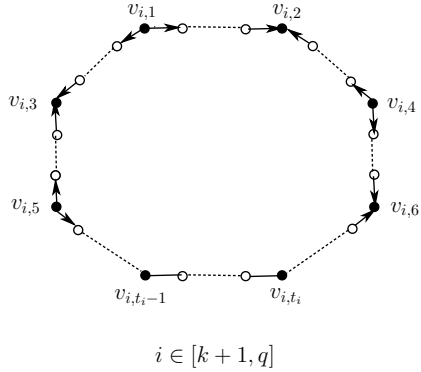}
\end{center}
\vskip-0.2cm
\caption{The orientation of an even cycle}\label{F2}
\end{figure}

Since $G_i$ is an $2d$-regular connected graph with $d \geq 2$ and $i\in[1,q]$, $G_{i}$ is an Euler graph. Let $C^{\ast}_{i}$ be an Euler tour of $G_{i}$. For each vertex $u\in V(G_{i})$, $C^{\ast}_{i}$ should pass through each vertex $d$ times. Let $u\in V(G_{i})$. Pick a fixed one of the $d$ copies of $u$ on $C^{\ast}_{i}$ as a real vertex and the remaining $d-1$ copies of $u$ as imaginary vertices. Then regarding $C^{\ast}_{i}$ as a circuit, say $C_i$, with $|G_{i}|$ real vertices. For any $i\in[1,q]$, we may assume that $|G_{i}|=t_{i}$. Let $n_{i}=\sum^{i}_{j=1}|E(C_{j})|$. Then $|E(C_{i})|=n_{i}-n_{i-1}$. Let $V^{i}_{R}=\{v_{i,1},v_{i,2},\ldots, v_{i,t_{i}}\}$ for $i\in[1,q]$. Then $V_{R}=\bigcup^{q}_{i=1}V^{i}_{R}$ is the set of real vertices of $\bigcup^{q}_{i=1}C_{i}$, where $v_{i,j}$ denote the $j$th real vertex of $C_{i}$. Let $V_{I}=V(\bigcup^{q}_{i=1}C_{i})\setminus V_{R}$ be the set of imaginary vertices of $\bigcup^{q}_{i=1}C_{i}$. By renaming the vertices in $V_{R}$ if necessary, we label the vertices of $V_{R}$ on $C_{i}$ with $v_{i,1},v_{i,2},v_{i,4},\ldots,v_{i,t_{i}-1},v_{i,t_{i}},v_{i,t_{i}-2},\ldots,v_{i,5},v_{i,3}$ in clockwise,
if $i\in[1,k]$ as depicted in Figure~\ref{F1}; and $v_{i,1},v_{i,2},v_{i,4},\ldots,v_{i,t_{i}},v_{i,t_{i}-1},v_{i,t_{i}-3},\ldots,v_{i,5},v_{i,3}$ in clockwise
if $i\in[k+1,q]$ as depicted in Figure~\ref{F2}.

Let $P^{i}_{j,k}(j<k)$ be the path between $v_{i,j}$ and $v_{i,k}$ on $C_{i}$ such that all internal vertices of $P^{i}_{j,k}$ are
not real vertices.
Next, we will find an orientation $D$ of $\bigcup^{q}_{i=1} C_{i}$.

When $i\in[1,k]$, set $d^{+}_{D}(v_{i,1})=1$, $d^{+}_{D}(v_{i,j})\in\{0,2\}$ for any $j\in\{2,3,\ldots,t_{i}\}$, and $d^{+}_{D}(u_{j})=1$ for each $u_{j}\in V_{I}\cap V(C_{i})$ by orienting the path $P^{i}_{1,2}$ from $v_{i,1}$ to $v_{i,2}$, the path $P^{i}_{2,4}$ from $v_{i,4}$ to $v_{i,2},\ldots$, $P^{i}_{3,5}$ from $v_{i,3}$ to $v_{i,5}$ and $P^{i}_{1,3}$ from $v_{i,3}$ to $v_{i,1}$ (see Figure~\ref{F1}).

When $i\in[k+1,q]$, set $d^{+}_{D}(v_{i,j})\in\{0,2\}$ for any $j\in[1,t_{i}]$, and $d^{+}_{D}(u_{j})=1$ for each $u_{j}\in V_{I}\cap V(C_{i})$ by orienting the path $P^{i}_{1,2}$ from $v_{i,1}$ to $v_{i,2}$, the path $P^{i}_{2,4}$ from $v_{i,4}$ to $v_{i,2}$, $\ldots$, $P^{i}_{3,5}$ from $v_{i,5}$ to $v_{i,3}$ and $P^{i}_{1,3}$ from $v_{i,1}$ to $v_{i,3}$ (see Figure~\ref{F2}).


 Assume that each edge on $C^{\ast}_{i}$ is oriented in the same as its oriention on $C_{i}$, then the orientation, say $D^{\ast}$, is an corresponding orientation of $G$. We need to find a labeling $c: A(D)\rightarrow[1,n_{q}]$ such that $c$ is an antimagic labeling of orientation $D$. Clearly, $c$ is also a labeling of $D^{\ast}$.

\begin{algorithm}
\caption{Label the edges of $C_{k+1},\ldots,C_{q}$}\label{alg1}
\begin{algorithmic}
\REQUIRE Even cycles $C_{i}$ for $i\in[k+1,q]$ with the given orientation $D$
\ENSURE A bijection $c_{e}:A(C_{k+1}\cup\ldots\cup C_{q})\rightarrow[n_{k}+1,n_{q}]$
\FOR{$i=k+1$ to $q$}
\STATE Assign the numbers in $[n_{i-1}+1,n_{i-1}+|P^{i}_{1,2}|]$ to the edges of $P^{i}_{1,2}$ in the increasing order along the orientation of $P^{i}_{1,2}$;\\
Assign the numbers in $[n_{i-1}+|P^{i}_{1,2}|+1,n_{i-1}+|P^{i}_{1,2}|+|P^{i}_{1,3}|]$ to the edges of $P^{i}_{1,3}$ in the increasing order along the orientation of $P^{i}_{1,3}$;\\
Set $V_{i}=\{v_{i,1},v_{i,2},v_{i,3}\}$;\\
\WHILE{$V_{i}\neq \{v_{i,1},v_{i,2},\ldots,v_{i,t_{i}}\}$}
\FOR{$j=2$ to $t_{i}-2$}
\STATE Assign the numbers in $[n_{i-1}+|P^{i}_{1,2}|+|P^{i}_{1,3}|+\ldots+|P^{i}_{j-1,j+1}|+1,n_{i-1}+|P^{i}_{1,2}|+|P^{i}_{1,3}|+\ldots+|P^{i}_{j,j+2}|]$ to the edges of $P^{i}_{j,j+2}$ in the increasing order along the orientation of $P^{i}_{j,j+2}$;\\
Set $V_{i}$ to be $V_{i}\cup\{v_{i,j+2}\}$;
\ENDFOR
\\Assign the numbers in $[n_{i-1}+|P^{i}_{1,2}|+|P^{i}_{1,3}|+\ldots+|P^{i}_{t_{i}-2,t_{i}}|+1,n_{i-1}+|P^{i}_{1,2}|+|P^{i}_{1,3}|+\ldots+|P^{i}_{t_{i}-2,t_{i}}|+|P^{i}_{t_{i}-1,t_{i}}|]$ to the edges of $P^{i}_{t_{i}-1,t_{i}}$ in the increasing order along the orientation of $P^{i}_{t_{i}-1,t_{i}}$;
\ENDWHILE
\ENDFOR
\end{algorithmic}
\end{algorithm}

The edges of even cycles $C_{i}$ are labeled as depicted in Algorithm 1 for $i\in[k+1,q]$.

For each $i\in[1,k]$ with $k\leq 9$, the bijection $c_{i}: A(C_{i})\rightarrow[n_{i-1}+1,n_{i}]$ and $c^{\prime}:A(C_{10}\cup\ldots\cup C_{k})\rightarrow[n_{9}+1,n_{k}]$ for $k\geq 10$ will be given
 such that $c$ is the desired labeling, where $c=c_{1}\cup c_{2}\cup\ldots\cup c_{k}\cup c_{e}$ if $k\leq 9$ and $c=c_{1}\cup c_{2}\cup\ldots\cup c_{9}\cup c^{\prime}\cup c_{e}$ for $k\geq 10$.

If $k=0$, we have $c=c_{e}$. If $k=1$, we define the bijection $c_{1}:A(C_{1})\rightarrow[1,n_{1}]$ as stated in Algorithm 2.

\begin{algorithm}
\caption{Label the edges of $C_{1}$}\label{alg2}
\begin{algorithmic}
\REQUIRE The odd cycle $C_{1}$ with the given orientation $D$
\ENSURE A bijection $c_{1}:A(C_{1})\rightarrow[1,n_{1}]$.\\
 Assign the numbers in $[1,|P^{1}_{1,2}|]$ to the edges of $P^{1}_{1,2}$ in the increasing order along the orientation of $P^{1}_{1,2}$;\\
 Assign the numbers in $[|P^{1}_{1,2}|+1,|P^{1}_{1,2}|+|P^{1}_{1,3}|]$ to the edges of $P^{1}_{1,3}$ in the increasing order along the orientation of $P^{1}_{1,3}$;\\
Set $V_{1}=\{v_{1,1},v_{1,2},v_{1,3}\}$;\\
\WHILE{$V_{1}\neq \{v_{1,1},v_{1,2},\ldots,v_{1,t_{1}}\}$}
\FOR{$j=2$ to $t_{1}-2$}
\STATE Assign the numbers in $[|P^{1}_{1,2}|+|P^{1}_{1,3}|+\ldots+|P^{1}_{j-1,j+1}|+1,|P^{1}_{1,2}|+|P^{1}_{1,3}|+\ldots+|P^{1}_{j,j+2}|]$ to the edges of $P^{1}_{j,j+2}$ in the increasing order along the orientation of $P^{1}_{j,j+2}$;\\
Set $V_{1}$ to be $V_{1}\cup\{v_{1,j+2}\}$;
\ENDFOR
\\Assign the numbers in $[|P^{1}_{1,2}|+|P^{1}_{1,3}|+\ldots+|P^{1}_{t_{1}-2,t_{1}}|+1,|P^{1}_{1,2}|+|P^{1}_{1,3}|+\ldots+|P^{1}_{t_{1}-2,t_{1}}|+|P^{1}_{t_{1}-1,t_{1}}|]$ to the edges of $P^{1}_{t_{1}-1,t_{1}}$ in the increasing order along the orientation of $P^{1}_{t_{1}-1,t_{1}}$;
\ENDWHILE
\end{algorithmic}
\end{algorithm}

By Algorithm 2, observe that the edges are labeled in $C_{1}$ in the order of $P^{1}_{1,2},P^{1}_{1,3},P^{1}_{2,4},P^{1}_{3,5},$ $P^{1}_{4,6},...,P^{1}_{t_{1}-2,t_{1}},P^{1}_{t_{1}-1,t_{1}}$ by using the numbers in $[1,n_{1}]$ with the increasingly order along the orientation of each path. If $k\in[2,9]$, we modify the label order of some fixed paths based on the Algorithm 2 to define the bijections $c_{2},...,c_{9}$. That is, when $k=2$, we label the edges in $C_{2}$ in the order of $P^{2}_{1,3},P^{2}_{1,2},P^{2}_{2,4},P^{2}_{3,5},P^{2}_{4,6},...,P^{2}_{t_{2}-2,t_{2}},P^{2}_{t_{2}-1,t_{2}}$ by using the numbers in $[n_{1}+1,n_{2}]$ with the increasingly order along the orientation of each path; when $3\leq k\leq 8$, for every $i\in[3,k]$, we label the edges in $C_{i}$ in the order of $P^{i}_{1,3},P^{i}_{2,4},P^{i}_{1,2},P^{i}_{3,5},P^{i}_{4,6},...,P^{i}_{t_{i}-2,t_{i}},P^{i}_{t_{i}-1,t_{i}}$ by using the numbers in $[n_{i-1}+1,n_{i}]$ with the increasingly order along the orientation of each path; when $k=9$, we label the edges in $C_{9}$ in the order of $P^{9}_{1,3},P^{9}_{2,4},P^{9}_{3,5},P^{9}_{1,2},P^{9}_{4,6},...,P^{9}_{t_{9}-2,t_{9}},P^{9}_{t_{9}-1,t_{9}}$ by using the numbers in $[n_{8}+1,n_{9}]$ with the increasing order along the orientation of each path. If $k\geq10$, we define the bijection $c^{\prime}$ such that $c^{\prime}:A(C_{10}\cup\ldots\cup C_{k})\rightarrow[n_{9}+1,n_{k}]$ is the same as stated in Algorithm 3.

\begin{algorithm}
\caption{For $k\geq10$, label the edges of $C_{10},\ldots,C_{k}$}\label{alg3}
\begin{algorithmic}
\REQUIRE Odd cycles $C_{i}$ for $i\in[10,k]$ with the given orientation $D$
\ENSURE A bijection $c^{\prime}:A(C_{10}\cup\ldots\cup C_{k})\rightarrow[n_{9}+1,n_{k}]$
\FOR{$i=10$ to $k$}
\STATE Assign the numbers in $[n_{i-1}+1,n_{i-1}+|P^{i}_{1,2}|]$ to the edges of $P^{i}_{1,2}$ in the increasing order along the orientation of $P^{i}_{1,2}$;\\
Assign the numbers in $[n_{i-1}+|P^{i}_{1,2}|+1,n_{i-1}+|P^{i}_{1,2}|+|P^{i}_{1,3}|]$ to the edges of $P^{i}_{1,3}$ in the increasing order along the orientation of $P^{i}_{1,3}$;\\
Set $V_{i}=\{v_{i,1},v_{i,2},v_{i,3}\}$;\\
\WHILE{$V_{i}\neq \{v_{i,1},v_{i,2},\ldots,v_{i,t_{i}}\}$}
\FOR{$j=2$ to $t_{i}-2$}
\STATE Assign the numbers in $[n_{i-1}+|P^{i}_{1,2}|+|P^{i}_{1,3}|+\ldots+|P^{i}_{j-1,j+1}|+1,n_{i-1}+|P^{i}_{1,2}|+|P^{i}_{1,3}|+\ldots+|P^{i}_{j,j+2}|]$ to the edges of $P^{i}_{j,j+2}$ in the increasing order along the orientation of $P^{i}_{j,j+2}$;\\
Set $V_{i}$ to be $V_{i}\cup\{v_{i,j+2}\}$;
\ENDFOR
\\Assign the numbers in $[n_{i-1}+|P^{i}_{1,2}|+|P^{i}_{1,3}|+\ldots+|P^{i}_{t_{i}-2,t_{i}}|+1,n_{i-1}+|P^{i}_{1,2}|+|P^{i}_{1,3}|+\ldots+|P^{i}_{t_{i}-2,t_{i}}|+|P^{i}_{t_{i}-1,t_{i}}|]$ to the edges of $P^{i}_{t_{i}-1,t_{i}}$ in the increasing order along the orientation of $P^{i}_{t_{i}-1,t_{i}}$;
\ENDWHILE
\ENDFOR
\end{algorithmic}
\end{algorithm}

It remains to verify that the bijection $c$ is an antimagic labeling of $D^{\ast}$, where  $c=c_{1}\cup c_{2}\cup\ldots\cup c_{k}\cup c_{e}$ if $k\leq 9$ and $c=c_{1}\cup c_{2}\cup\ldots\cup c_{9}\cup c^{\prime}\cup c_{e}$ for $k\geq 10$.

Since $C_{i}$ corresponds to $C^{\ast}_{i}$, and $C^{\ast}_{i}$ can be reselected if necessary, so $C_{i}$ can satisfy the following conditions according to the different values of $k$, respectively.

$(1)$ For $k\in[0,5d+4]$, if $k\in[3,6]$, let $|P^{i}_{2,4}|=i-2$ for $i\in[3,k]$;
If $k=7$ or $k=8$, let $|P^{i}_{2,4}|=i-2$ for $i\in[3,k]$ and $|P^{1}_{2,4}|\geq 3$;
If $k=9$, let $|P^{i}_{2,4}|=i-2$ for $i\in[3,9]$, $|P^{9}_{3,5}|=1$ and $|P^{1}_{2,4}|\geq 4$;
If $k\in[10,5d+4]$, let $|P^{1}_{1,2}|=|P^{1}_{1,3}|=|P^{9}_{3,5}|=|P^{i}_{1,3}|=1$ for $i\in[10,k]$, $|P^{i}_{2,4}|=i-2$ for $i\in[3,9]$, $|P^{i}_{1,2}|=i-8$ for $i\in[10,k]$ and $|P^{1}_{2,4}|\geq 5d-6$.

$(2)$ For $k\geq5d+5$, based on the orientation and $|G_{1}|\geq 2x_{0}+5$, we can let $|P^{1}_{1,2}|=|P^{1}_{1,3}|=|P^{9}_{3,5}|=|P^{i}_{1,3}|=1$ for $i\in[10,k]$, $|P^{i}_{2,4}|=i-2$ for $i\in[3,9]$, $|P^{i}_{1,2}|=i-8$ for $i\in[10,k]$ and $|P^{1}_{2,4}|\geq (2d-2)x_{0}+5d-6$.\\

The definition $V_{R}$ and $V_{I}$ and the labeling method of even cycles are the same as the method in~\cite{L}.

\begin{clm}\label{clm1}
If $D$ is antimagic, then $D^{\ast}$ is antimagic.
\end{clm}
\f {\bf Proof of Claim 1.} By the three algorithms, $s_{D}(u_{j})=-1$ for all $u_{j}\in V_{I}$. We may assume that $V(G)=V_{R}$. For each $v\in V(G),s_{D^{\ast}}(v)=s_{D}(v)+(d-1)s_{D}(u^{\ast})=s_{D}(v)-(d-1)$, where $u^{\ast}$ is one of the $d-1$ imaginary vertices of $v$. Therefore, if $D$ is antimagic, for any $u,v\in V_{R}$ with $u\neq v$, $s_{D}(u)\neq s_{D}(v)$. Then, for any $u,v\in V(G)$ with $u\neq v$, one has that $s_{D^{\ast}}(u)\neq s_{D^{\ast}}(v)$. That is, $D^{\ast}$ is antimagic.
\hfill\qed

By Claim 1, it suffices to show that for any $u,v\in V_{R}$ with $u\neq v$, $s_{D}(u)\neq s_{D}(v)$.

\begin{clm}\label{clm2}
For any $u$ and $v\in V(C_{k+1}\cup \ldots\cup C_{q})\cap V_{R}$ with $u\neq v$, one has that $s_{D}(u)\neq s_{D}(v)$.
\end{clm}
\f {\bf Proof of Claim 2.} By the orientation of even cycle, observe that the two edges incident with each real vertex be either both entering the vertex or both leaving the vertex. Choose two different real vertices $u$ and $v$ from all even cycles. Clearly, if $d^{+}_{D}(u)=0$ and $d^{+}_{D}(v)=2$ or $d^{+}_{D}(v)=0$ and $d^{+}_{D}(u)=2$, we have $s_{D}(u)\cdot s_{D}(v)<0$. If $d^{+}_{D}(u)=d^{+}_{D}(v)=0$ or $d^{+}_{D}(u)=d^{+}_{D}(v)=2$, by Algorithm 1, the labels of two edges incident with one vertex must be strictly less than the labels of two edges incident with the other vertex, respectively. One has that $|s_{D}(u)|<|s_{D}(v)|$ or $|s_{D}(v)|<|s_{D}(u)|$. Thus, $s_{D}(u)\neq s_{D}(v)$.
\hfill\qed

\begin{clm}\label{clm3}
For any $u\in V(C_{1}\cup\ldots\cup C_{k})\cap V_{R}$ and $v\in V(C_{k+1}\cup\ldots\cup C_{q})\cap V_{R}$, one has that $s_{D}(u)\neq s_{D}(v)$.
\end{clm}
\f {\bf Proof of Claim 3.} By the definition of the labeling $c$, $A(C_{1}\cup\ldots\cup C_{k})\rightarrow[1,n_{k}]$ and $A(C_{k+1}\cup\ldots\cup C_{q})\rightarrow[n_{k}+1,n_{q}]$. The labels of two edges incident with $u$ must be strictly less than the labels of two edges incident with $v$, respectively. If $u$ is in an odd cycle $C_i$ but $u$ is not the first vertex $v_{i,1}$ of $C_{i}$, it is  clearly that $s_{D}(u)\neq s_{D}(v)$. If $u=v_{i,1}$ in $C_{i}$ for $i\in[1,k]$, assume that the labels of the edges entering $u$ is $a$ and leaving $u$ is $b$. Assume that the labels of two edges incident with $v$ are $c$ and $d$ respectively. If $s_{D}(u)\cdot s_{D}(v)<0$, then $s_{D}(u)\neq s_{D}(v)$. If $s_{D}(u)=a-b>0$ and $s_{D}(v)=c+d>0$, then $a-b<a+b<c+d$, so $s_{D}(u)\neq s_{D}(v)$. Otherwise, $s_{D}(u)=a-b<0$ and $s_{D}(v)=-(c+d)<0$, one has that $-(c+d)<-(a+b)<a-b$, it also implies that $s_{D}(u)\neq s_{D}(v)$.
\hfill\qed

By Claims 2 and 3, it suffices to prove Claim~\ref{clm4}.

\begin{clm}\label{clm4}
For any $u$ and $v\in V(C_{1}\cup \ldots\cup C_{k})\cap V_{R}$ with $u\neq v$, one has that $s_{D}(u)\neq s_{D}(v)$.
\end{clm}
\f {\bf Proof of Claim 4.}
The proof will be given by the induction on the number of odd cycles $k\geq 2$.

If $k=2$, for any $u$ and $v\in V(C_{1}\cup C_{2})\cap V_{R}$, one has that $s_{D}(u)\neq s_{D}(v)$, which has been proved in~\cite{L}, as the labeling method of $V(C_{1}\cup C_{2})$ is the same as the method in~\cite{L}.

Assume that the number of odd cycles is not more than $k-1$, the result is true. Now suppose that there are $k$ odd cycles with $k\geq 3$.

 If $u$ and $v\in V(C_{1}\cup\ldots\cup C_{k-1})\cap V_{R}$, since $D\setminus C_k$ contains $k-1$ odd cycles, so $s_{D\setminus C_k}(u)\neq s_{D\setminus C_k}(v)$ by the inductive hypothesis. Since $s_{D\setminus C_k}(u)=s_{D}(u)$ and $s_{D\setminus C_k}(v)=s_{D}(v)$, so $s_{D}(u)\neq s_{D}(v)$. Thus we need to consider only the following two cases.

 Case 1. $u,v\in V(C_{k})\cap V_{R}$. Without loss of generality, let $u=v_{k,t}$ and $v=v_{k,\ell}$ and $t<\ell$.

 If $s_{D}(u)\cdot s_{D}(v)<0$, one has that $s_{D}(u)\neq s_{D}(v)$. So we only need to consider $s_{D}(u)\cdot s_{D}(v)>0$.

 By the labeling of $C_{k}$, the labels of two edges incident with $u$ must be strictly less than the labels of two edges incident with $v$, respectively. If $s_{D}(u)$ and $s_{D}(v)$ are both positive, then $s_{D}(u)<s_{D}(v)$; if $s_{D}(u)$ and $s_{D}(v)$ are both negative, then $s_{D}(u)>s_{D}(v)$. Thus, $s_{D}(u)\neq s_{D}(v)$.

 Case 2. $u\in V(C_{1}\cup\ldots\cup C_{k-1})\cap V_{R}$ and $v\in V(C_{k})\cap V_{R}$.

 Let $X=\{v_{1,1},v_{2,1},...,v_{k-1,1}\}$. There are the following three subcases.

 Subcase 2.1. $u\in X$ and $v=v_{k,1}$.

 Since $s_{D}(v_{1,1})=1,s_{D}(v_{j,1})=-(j-1)$ for $j\in[2,9]$, and $s_{D}(v_{j,1})=j-8$ for $j\in[10,k]$, we have $s_{D}(u)\neq s_{D}(v)$.

 Subcase 2.2. $u\in V(C_{1}\cup\ldots\cup C_{k-1})\cap V_{R}$ and $v\in V(C_{k})\cap V_{R}\setminus\{v_{k,1}\}$.

 By the definition of $c$, $A(C_{1}\cup\cdots\cup C_{k-1})\rightarrow [1,n_{k-1}]$ and $A(C_{k})\rightarrow [n_{k-1}+1,n_{k}]$. Using the similar arguments in the proof of Claim 3, we can easily show that $s_{D}(u)\neq s_{D}(v)$.

Subcase 2.3. $u\in V(C_{1}\cup\ldots\cup C_{k-1})\cap V_{R}\setminus X$ and $v=v_{k,1}$.

If $k\in[3,9]$, then $s_{D}(v_{k,1})=-(k-1)<0$. If $s_{D}(u)>0$, we have $s_{D}(u)\neq s_{D}(v)$. Let us consider the set, say $W_{1}$, of all vertices $w$ in $V(C_{1}\cup\ldots\cup C_{k-1})\cap V_{R}\setminus X$ such that $s_{D}(w)<0$. Note that among all the vertices of $W_{1}$, the labels of two edges incident with $v_{1,3}$ are minimum. By the conditions, one has that
 $|s_{D}(v_{1,3})|=(|P^{1}_{1,2}|+1)+(|P^{1}_{1,2}|+|P^{1}_{1,3}|+|P^{1}_{2,4}|+1)=|P^{1}_{2,4}|+5\geq 9>k-1=s_{D}(v_{k,1})$. That is, $s_{D}(w)> s_{D}(v_{1,3})>s_{D}(v_{k,1})$ for all $w\in W_{1}$. Thus, if $s_{D}(u)<0$, one still has that $s_{D}(u)\neq s_{D}(v)$.

If $k\in[10,5d+4]$, then $s_{D}(v_{k,1})=k-8>0$. If $s_{D}(u)<0$, one has that $s_{D}(u)\neq s_{D}(v)$.
Let us consider the set, say $W_{2}$, of all vertices $w$ in $V(C_{1}\cup\ldots\cup C_{k-1})\cap V_{R}\setminus X$ such that $s_{D}(w)>0$.
Note that among all the vertices of $W_{2}$, the labels of two edges incident with $v_{1,2}$ are minimum. By the conditions that $|P^{1}_{2,4}|\geq 5d-6$, one has that $s_{D}(v_{1,2})=|P^{1}_{1,2}|+(|P^{1}_{1,2}|+|P^{1}_{1,3}|+|P^{1}_{2,4}|)=|P^{1}_{2,4}|+3\geq 5d-3>k-8=s_{D}(v_{k,1})$. That is, $s_{D}(v_{k,1})<s_{D}(v_{1,2})< s_{D}(w)$ for all $w\in W_{2}$. Thus, if $s_{D}(u)>0$, we still have $s_{D}(u)\neq s_{D}(v)$.

If $k\geq 5d+5$, then $s_{D}(v_{k,1})=k-8>0$. If $s_{D}(u)<0$, we have $s_{D}(u)\neq s_{D}(v)$.
Let us consider the set, say $W_{3}$, of all vertices $w$ in $V(C_{1}\cup\ldots\cup C_{k-1})\cap V_{R}\setminus X$ such that $s_{D}(w)>0$. Note that among all the vertices of $W_{3}$, the labels of two edges incident with $v_{1,2}$ are minimum.
If $k\geq 5d+5$, by the condition of Theorem~\ref{thm1}, $|G_{1}|\geq 2x_{0}+5$, where $x_{0}$ is a unique positive integer solution for one of equations: $k=(2d-2)(x+2)+0$, $k=(2d-2)(x+2)+1$,\ldots, $k=(2d-2)(x+2)+(d+8), k=(2d-2)(x+1)+(d+9)$, $k=(2d-2)(x+1)+(d+10),\ldots, k=(2d-2)(x+1)+(2d-3)$. It implies that, $x_{0}$ is a positive integer solution for one of equations: $k=(2d-2)x+3d+7,k=(2d-2)x+3d+8, \ldots,k=(2d-2)x+5d+4$. So, $k\leq(2d-2)x_{0}+5d+4$. (Otherwise, $k>(2d-2)x_{0}+5d+4$. Assume $k=(2d-2)x_{0}+5d+t$ for $t\geq 5$.
It implies that $k=(2d-2)(x_{0}+1)+3d+t+2$ with $3d+t+2\geq 3d+7$.
It implies that $x_{0}+1$ or $x_{0}+i$ for $i>1$ is the positive integer solution, which contradicts with $x_{0}$ is the unique positive integer solution.) Since $s_{D}(v_{1,2})=|P^{1}_{1,2}|+(|P^{1}_{1,2}|+|P^{1}_{1,3}|+|P^{1}_{2,4}|)\geq |P^{1}_{2,4}|+3$ and $|P^{1}_{2,4}|\geq (2d-2)x_{0}+5d-6$ given in the definition of $c$, $s_{D}(v_{1,2})\geq (2d-2)x_{0}+5d-3>k-8=s_{D}(v_{k,1})$. That is, $s_{D}(v_{k,1})<s_{D}(v_{1,2})<s_{D}(w)$ for all $w\in W_{3}$. Thus, if $s_{D}(u)>0$, we still have $s_{D}(u)\neq s_{D}(v)$.

Therefore, the result holds if the number of odd cycles is $k$ with $k\geq3$.

This completes the proof of Theorem~\ref{thm1}.
\hfill\qed

\section*{Acknowledgments}
This work was supported by the National Natural Science Foundation of China
(No. 11731002), the Fundamental Research Funds for the Central Universities (No. 2016JBM071, 2016JBZ012).


\begin{thebibliography}{6}

\bibitem{A} N. Alon, G. Kaplan, A. Lev, T.Roditty, and R. Yuster, Dense graphs are antimagic, J. Graph Theory 47(2004) 297-309.

\bibitem{CH} F. Chang, Y-C. Liang, Z. Pan, and X. Zhu, Antimagic labeling of regular graphs, J. Graph Theory 82(2016) 339-349.

\bibitem{CR} D.W. Cranston, Regular bipartite graphs are antimagic, J Graph Theory 60(2009) 173-182.

\bibitem{CRL} D.W. Cranston, Y-C. Liang, and X. Zhu, Regular graphs of odd degree are antimagic, J. Graph Theory 80(2015) 28-33.

\bibitem{H} D. Hefetz, T. M$\ddot{u}$tze, and J. Schwartz, On antimagic directed graphs, J. Graph Theory 64(2010) 219-232.

\bibitem{L} T. Li, Z.-X. Song, G.H. Wang, D.L. Yang and C.-Q. Zhang, Antimagic orientations of even regular graphs, Available at arXiv:1707.03507.

\bibitem{S} S. Shan and X. Yu, Antimagic orientation of biregular bipartite graphs, The Electronic Journal of Combinatorics 24(4)(2017).

\end{thebibliography}
\end{document}